\documentclass[11pt]{amsart}

\author[C.~Sanna]{Carlo Sanna$^\dagger$}
\thanks{$\dagger\,$C.~Sanna is a member of GNSAGA of INdAM and of CrypTO, the group of Cryptography and Number~Theory of Politecnico di Torino}
\address{\parbox{\linewidth}{
Politecnico di Torino, Department of Mathematical Sciences\\
Corso Duca degli Abruzzi 24, 10129 Torino, Italy\\[-8pt]}}
\email{carlo.sanna.dev@gmail.com}


\title{Membership in random ratio sets}

\usepackage{amsmath}
\usepackage{amssymb}
\usepackage{amsthm}
\usepackage{geometry}
\usepackage{placeins}
\usepackage[colorlinks=true]{hyperref}
\usepackage{graphicx}
\usepackage{bm}

\newtheorem{theorem}{Theorem}[section]
\newtheorem{corollary}{Corollary}[section]

\theoremstyle{remark}
\newtheorem{remark}{Remark}[section]
\newtheorem{example}{Example}[section]

\DeclareMathOperator{\lcm}{lcm}
\DeclareMathOperator{\Li}{Li}
\DeclareMathOperator{\vis}{vis}

\newcommand{\Beta}{\text{B}}

\begin{document}

\begin{abstract}
Let $\mathcal{A}$ be a random set constructed by picking independently each element of $\{1, \dots, n\}$ with probability $\alpha \in (0, 1)$.
We give a formula for the probability that a rational number $q$ belong to the random ratio set $\mathcal{A} /\! \mathcal{A} := \{a / b : a,b \in \mathcal{A}\}$.
This generalizes a previous result of Cilleruelo and Guijarro-Ord\'{o}\~{n}ez.
Moreover, we make some considerations about formulas for the probability of the event $\bigvee_{i=1}^k\!\big(q_i \in \mathcal{A} /\! \mathcal{A}\big)$, where $q_1, \dots, q_k$ are rational numbers, showing that they are related to the study of the connected components of certain graphs.
In particular, we give formulas for the probability that $q^e \in \mathcal{A} /\! \mathcal{A}$ for some $e \in \mathcal{E}$, where $\mathcal{E}$ is a finite or cofinite set of positive integers with $1 \in \mathcal{E}$.
\end{abstract}

\maketitle

\section{Introduction}

For every positive integer $n$ and for every $\alpha \in (0, 1)$, let $\mathcal{B}(n, \alpha)$ denote the probabilistic model in which a random set $\mathcal{A} \subseteq \{1, \dots, n\}$ is constructed by picking independently every element of $\{1, \dots, n\}$ with probability $\alpha$.
Several authors studied number-theoretic objects involving random sets in this probabilistic model, including: the least common multiple $\lcm(\mathcal{A})$~\cite{MR4009436, MR3239153} (see also~\cite{MR4091939}), the product set $\mathcal{A}\mathcal{A} := \{ab : a, b \in \mathcal{A}\}$~\cite{MR3640773, MR4221525, MR4169019}, and the ratio set $\mathcal{A} /\!\mathcal{A} := \{a / b : a, b \in \mathcal{A}\}$~\cite{MR3649012, MR3640773}.

Regarding random ratio sets, Cilleruelo and Guijarro-Ord\'{o}\~{n}ez~\cite{MR3649012} proved the following:

\begin{theorem}\label{thm:cardinality}
Let $\mathcal{A}$ be a random set in $\mathcal{B}(n, \alpha)$.
Then, for $\alpha$ fixed and $n \to +\infty$, we have
\begin{equation*}
|\mathcal{A} /\!\mathcal{A}| \sim \frac{6}{\pi^2} \cdot \frac{\alpha^2 \Li_2(1 - \alpha^2)}{1 - \alpha^2} \cdot n^2 ,
\end{equation*}
with probability $1 - o(1)$, where $\Li_2(z) := \sum_{k = 1}^\infty z^k / k^2$ is the dilogarithm function.
\end{theorem}

A fundamental step in the proof of Theorem~\ref{thm:cardinality} is determining a formula for the probability that certain rational numbers belong to $\mathcal{A} /\!\mathcal{A}$.
Precisely, Cilleruelo and Guijarro-Ord\'{o}\~{n}ez~\cite[Eq.~(2)]{MR3649012} showed that for all positive integers $r < s$, with $(r, s) = 1$ and $s > n^{1/2}$, we have 
\begin{equation*}
\mathbb{P}\big(r / s \in \mathcal{A} /\!\mathcal{A}\big) = 1 - (1 - \alpha^2)^{\lfloor n / s \rfloor} .
\end{equation*}
Note that the assumption $r < s$ is not restrictive, since $r / s \in \mathcal{A} /\! \mathcal{A}$ if and only if $s / r \in \mathcal{A} /\! \mathcal{A}$, while the assumption $s > n^{1/2}$ is indeed a restriction.

Our first result is a general formula for the probability that a rational number belongs to the ratio set $\mathcal{A} /\! \mathcal{A}$.

\begin{theorem}\label{thm:prob}
Let $\mathcal{A}$ be a random set in $\mathcal{B}(n, \alpha)$.
Then we have
\begin{equation}\label{equ:prob}
\mathbb{P}\big(r / s \in \mathcal{A} /\! \mathcal{A}\big) = 1 - \prod_{i \,=\, 1}^{\left\lfloor \frac{\log n}{\log s}\right\rfloor} \gamma_i^{\lfloor n / s^i \rfloor} ,
\end{equation}
for all positive integers $r < s$ with $(r, s) = 1$, where $\gamma_i := \beta_{i - 1} \beta_{i + 1} / \beta_i^2$ with $\beta_0 := 1$, $\beta_1 := 1$, and $\beta_{i + 1} := (1 - \alpha)\beta_i + \alpha (1 - \alpha) \beta_{i-1}$, for all integers $i \geq 1$.
\end{theorem}

\begin{remark}
If $\alpha = 1 / 2$ then for all integers $i \geq 0$ we have $\beta_i = F_{i + 2} / 2^i$, where $\{F_i\}_{i=0}^\infty$ is the sequence of Fibonacci numbers, defined recursively by $F_0 := 0$, $F_1 := 1$, and $F_{i + 2} := F_{i + 1} + F_i$. 
\end{remark}

As a consequence of Theorem~\ref{thm:prob}, we obtain the following corollary:

\begin{corollary}\label{cor:prob}
Let $\mathcal{A}$ be a random set in $\mathcal{B}(n, \alpha)$.
Then we have
\begin{equation}
\mathbb{P}\big(r / s \in \mathcal{A} /\! \mathcal{A}\big) = 1 - \exp\!\big({-}\delta(s)n + O_\alpha(1)\big) ,
\end{equation}
for all positive integers $r < s \leq n$ with $(r, s) = 1$, where
\begin{equation}\label{equ:delta-series}
\delta(s) := \sum_{i \,=\, 1}^\infty \frac{\log (1/\gamma_i)}{s^i}
\end{equation}
is an absolutely convergent series.
\end{corollary}

It is natural to ask if Theorem~\ref{thm:prob} can be generalized to a formula for the probability of the event $\bigvee_{i=1}^k \big(r_i / s_i \in \mathcal{A} /\! \mathcal{A}\big)$, where $r_1/s_1, \dots, r_k / s_k$ are rational numbers.
The answer should be ``yes'', but the task seems very complex (see Section~\ref{sec:generalizations} for more details).

However, we proved the following result concerning powers of the same rational number.

\begin{theorem}\label{thm:powers}
Let $\mathcal{A}$ be a random set in $\mathcal{B}(n, \alpha)$, and let $\mathcal{E}$ be a finite or cofinite set of positive integers with $1 \in \mathcal{E}$.
Then we have
\begin{equation}\label{equ:prob-powers}
\mathbb{P}\!\left(\,\bigvee_{e \,\in\, \mathcal{E}} \big((r / s)^e \in \mathcal{A} /\! \mathcal{A}\big)\right) = 1 - \prod_{i \,=\, 1}^{\left\lfloor \frac{\log n}{\log s}\right\rfloor} \left(\gamma_i^{(\mathcal{E})}\right)^{\lfloor n / s^i \rfloor} ,
\end{equation}
for all positive integers $r < s$ with $(r, s) = 1$, where $\gamma_i^{(\mathcal{E})} := \beta_{i - 1}^{(\mathcal{E})} \beta_{i + 1}^{(\mathcal{E})} / \big(\beta_i^{(\mathcal{E})}\big)^2$, for all integers $i \geq 1$, and $\big\{\beta_j^{(\mathcal{E})}\big\}_{j = 0}^\infty$ is a linear recurrence depending only on $\mathcal{E}$ and $\alpha$.
In particular, if $\mathcal{E}$ is cofinite then $\gamma_i^{(\mathcal{E})}$ is a rational function of $i$, for all sufficiently large $i$.
\end{theorem}

As a matter of example, we provide the following:

\begin{example}
$\beta_0^{\{1,2\}} = 1$, $\beta_1^{\{1,2\}} = 1$, $\beta_2^{\{1,2\}} = 1 - \alpha^2$, and 
\begin{equation*}
\beta_i^{\{1,2\}} = (1 - \alpha) \beta_{i - 1}^{\{1,2\}} + \alpha (1 - \alpha)^2 \beta_{i - 3}^{\{1,2\}}
\end{equation*}
for all integers $i \geq 3$.
\end{example}

\begin{example}
$\beta_0^{\{1,3\}} = 1$, $\beta_1^{\{1,3\}} = 1$, $\beta_2^{\{1,3\}} = 1 - \alpha^2$, $\beta_3^{\{1,3\}} = 1 - 2\alpha^2 + \alpha^3$, and
\begin{equation*}
\beta_i^{\{1,3\}} = (1 - \alpha) \beta_{i-1}^{\{1,3\}} + \alpha(1 - \alpha)\beta_{i-2}^{\{1,3\}} - \alpha(1 - \alpha)^2 \beta_{i-3}^{\{1,3\}} + \alpha(1 - \alpha)^3 \beta_{i-4}^{\{1,3\}}
\end{equation*}
for all integers $i \geq 4$.
\end{example}

\begin{example}
$\beta_i^{(\mathbb{N})} = (1 - \alpha)^{i-1} \big((1 - \alpha) + i \alpha\big)$ and
\begin{equation*}
\gamma_i^{(\mathbb{N})} = 1 - \frac1{(i + \alpha^{-1} - 1)^{2}}
\end{equation*}
for all integers $i \geq 1$.
\end{example}

\begin{example}
$\beta_i^{\mathbb{N} \setminus \{2\}} = (1 - \alpha)^{i-2}\big((1 - \alpha)^2 + i \alpha(1 - \alpha) + (i - 2)\alpha^2\big)$ and
\begin{equation*}
\gamma_i^{\mathbb{N} \setminus \{2\}} = 1 - \frac1{(i + \alpha^{-1} - \alpha - 2)^2}
\end{equation*}
for all integers $i \geq 3$.
\end{example}

\section{Notation}

We employ the Landau--Bachmann ``Big Oh'' notation $O$ with its usual meaning.
Any dependence of implied constants is explicitly stated or indicated with subscripts.
Greek letters are reserved for quantities that depends on $\alpha$.

\section{Proof of Theorem~\ref{thm:prob}}

\begin{figure}[h]
\includegraphics[width=\textwidth]{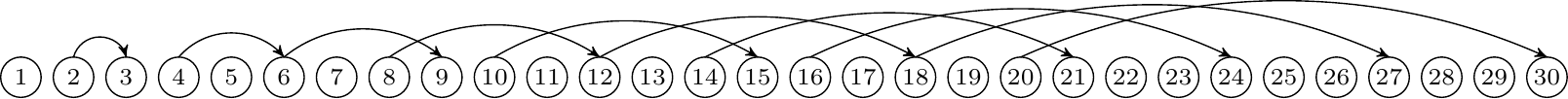}
\caption{The directed graph $\mathcal{G}(30;2,3)$.}\label{fig:G-30-2-3}
\end{figure}

Let us consider the directed graph $\mathcal{G}(n; r, s)$ having vertices $1, \dots, n$ and edges $rt \to st$, for all positive integers $t \leq n / s$.
For an example, see Figure~\ref{fig:G-30-2-3}.
Note that two edges $rt \to st$ and $rt^\prime \to st^\prime$, with $t < t^\prime$, have a common vertex if and only if $st = rt^\prime$.
In such a case, recalling that $r$ and $s$ are relatively prime, it follows that $t = ru$ and $t^\prime = su$, for some positive integer~$u$, and consequently the two edges form a directed path $r^2 u \to rs u \to s^2 u$.
Indeed, iterating this reasoning, it follows that all paths of $i + 1$ vertices are of the type
\begin{equation*}
r^i u \to r^{i-1} s u \to \cdots \to r s^{i-1} u \to s^i u ,
\end{equation*}
for some positive integer $u \leq n / s^i$.
Morever, it is easy to check that each vertex of $\mathcal{G}(n; r, s)$ is incident to at most two edges.
Therefore, all the connected components of $\mathcal{G}(n; r, s)$ are directed path graphs of at most $k := \lfloor \log n /\! \log s\rfloor + 1$ vertices.

Let $c_i$, respectively $d_i$, be the number of connected components, respectively directed paths, of $\mathcal{G}(n;r,s)$ consisting of exactly $i$ vertices (considering each isolated vertex as a directed path with one vertex).
On the one hand, from the previous reasonings, we have that $d_i = \lfloor n / s^{i - 1} \rfloor$ for each positive integer $i$.
On the other hand, since each connected component of $i$ vertices contains exactly $i - j + 1$ directed paths of $j$ vertices, for all positive integers $j \leq i$, we have
\begin{equation}\label{equ:linear}
\sum_{i \,=\, j}^k (i - j + 1)\, c_i = d_j , \quad j = 1, \dots, k .
\end{equation}
The linear system~\eqref{equ:linear} in unknowns $c_1, \dots, c_k$ can be solved by subtracting to each equation the next one, and then again subtracting to each equation the next one.
This yields
\begin{equation}\label{equ:ci}
c_i = d_i - 2d_{i + 1} + d_{i + 2}, \quad i=1,\dots,k .
\end{equation}
(Note that $d_i = 0$ for every integer $i > k$.)

Now we have that $r / s \in \mathcal{A} /\! \mathcal{A}$ if and only if there exists a positive integer $t \leq n / s$ such that $rt \in \mathcal{A}$ and $st \in \mathcal{A}$.
Therefore,
\begin{equation*}
\mathbb{P}\big(r / s \notin \mathcal{A} /\! \mathcal{A}\big) = \mathbb{P}\!\left(\bigwedge_{t \,\leq\, n / s} E(rt, st)\right)
\end{equation*}
where $E(a,b)$ denotes the event $\lnot (a \in \mathcal{A} \land b \in \mathcal{A})$.
Clearly, there is a natural correspondence between the events $E(rt, st)$ and the edges $rt \to st$ of $\mathcal{G}(n; r, s)$.
In particular, events corresponding to edges of different connected components of $\mathcal{G}(n; r, s)$ are independent.
Furthermore, if $m_1 \to \cdots \to m_i$ is a connected component of $\mathcal{G}(n; r, s)$, then 
$$\mathbb{P}\big(E(m_1, m_2) \land \cdots \land E(m_{i-1}, m_i)\big)$$
is equal to the probability that the string $m_1, \dots, m_i$ has no consecutive elements in $\mathcal{A}$.
In~turn, this is easily seen to be equal to $\beta_i$.
Indeed, for $i=0,1$ the claim is obvious, since $\beta_0 = 1$ and $\beta_1 = 1$; while for $i \geq 2$ we have that $m_1, \dots, m_i$ contains no consecutive elements in $\mathcal{A}$ if and only if either $m_i \notin \mathcal{A}$ and $m_1, \dots, m_{i-1}$ has no consecutive elements in $\mathcal{A}$, or $m_i \in \mathcal{A}$, $m_{i - 1} \notin \mathcal{A}$, and $m_1, \dots, m_{i-2}$ has no consecutive elements in $\mathcal{A}$, so that the claim follows from the recursion $\beta_i = (1 - \alpha) \beta_{i - 1} + \alpha(1 - \alpha)\beta_{i - 2}$.

Therefore, also using~\eqref{equ:ci}, we have
\begin{align}\label{equ:chain}
\mathbb{P}\big(r / s \notin \mathcal{A} /\! \mathcal{A}\big) &= \prod_{i \,=\, 1}^k \prod_{\substack{m_1 \,\to\, \cdots \,\to\, m_i \\ \text{\,con.~com.~of~$\mathcal{G}(n;r,s)$}}} \mathbb{P}\big(E(m_1, m_2) \land \cdots \land E(m_{i-1}, m_i)\big) \\
&= \prod_{i \,=\, 1}^k \beta_i^{c_i} = \prod_{i \,=\, 1}^k \beta_i^{d_i - 2d_{i+1} + d_{i + 2}} = \prod_{i \,=\, 1}^{k-1} \left(\frac{\beta_{i-1}\beta_{i+1}}{\beta_i^2}\right)^{d_{i+1}} = \prod_{i \,=\, 1}^{\left\lfloor \frac{\log n}{\log s} \right\rfloor} \gamma_i^{\lfloor n / s^i \rfloor} , \nonumber
\end{align}
and~\eqref{equ:prob} follows.
The proof is complete.

\section{Proof of Corollary~\ref{cor:prob}}

Throughout this section, implied constants may depend on $\alpha$.
Let $\rho_1, \rho_2$ be the roots of the characteristic polynomial $X^2 - (1 - \alpha)X - \alpha(1-\alpha)$ of the linear recurrence $\{\beta_i\}_{i=0}^\infty$.
Recalling that $\alpha \in (0, 1)$, an easy computation shows that $|\rho_1| \neq |\rho_2|$.
Without loss of generality, assume $|\rho_1| > |\rho_2|$ and put $\varrho := |\rho_2 / \rho_1|$, so that $\varrho \in (0, 1)$.
Hence, there exist complex numbers $\zeta_1, \zeta_2$ such that
\begin{equation*}
\beta_i = \zeta_1 \rho_1^i + \zeta_2 \rho_2^i = \zeta_1 \rho_1^i \big(1 + O(\varrho^i)\big) ,
\end{equation*}
for every integer $i \geq 0$.
Consequently, we have
\begin{equation*}
\gamma_i = \frac{\beta_{i - 1}\beta_{i+1}}{\beta_i^2} = \frac{\zeta_1 \rho_1^{i-1} \big(1 + O(\varrho^{i-1})\big)\zeta_1 \rho_1^{i+1} \big(1 + O(\varrho^{i+1})\big)}{\left(\zeta_1 \rho_1^i\big(1 + O(\varrho^i)\big)\right)^2} = 1 + O(\varrho^i) ,
\end{equation*}
and $\log \gamma_i = O(\varrho^i)$, for every sufficiently large integer $i$.
In particular, it follows that~\eqref{equ:delta-series} is an absolutely convergent series.

Now put $\ell := \lfloor \log n /\! \log s \rfloor$.
From Theorem~\ref{thm:prob}, we get that
\begin{equation*}
\mathbb{P}\big(r / s \in \mathcal{A} /\! \mathcal{A}\big) = 1 - e^L ,
\end{equation*}
where
\begin{align*}
L &:= \sum_{i \,=\, 1}^\ell \left\lfloor \frac{n}{s^i} \right\rfloor\! \log \gamma_i = \sum_{i \,=\, 1}^\infty \frac{\log \gamma_i}{s^i}\, n + O\!\left(\,\sum_{i \,>\, \ell} \frac{|\!\log \gamma_i|}{s^i}\, n\right) + O\!\left(\,\sum_{i \,=\, 1}^\ell |\!\log \gamma_i| \right) \\
 &= {-}\delta(s)n + O\!\left(\frac{n}{s^{\ell+1}}\right) + O\!\left(\,\sum_{i \,=\, 1}^\ell \varrho^i \right) = {-}\delta(s)n + O(1) ,
\end{align*}
as desired.
The proof is complete.

\begin{remark}\label{rmk:alternating}
A more detailed analysis shows that $\left\{\gamma_i^{(-1)^i}\right\}_{i=1}^\infty$ is a strictly decreasing sequence tending to $1$. 
In particular,~\eqref{equ:delta-series} is an alternating series.
\end{remark}

\section{Proof of Theorem~\ref{thm:powers}}

\begin{figure}[h]
\includegraphics[width=\textwidth]{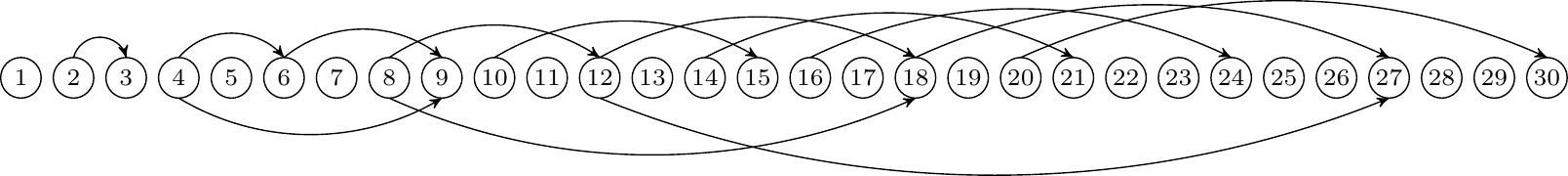}
\caption{The directed graph $\mathcal{G}^{\{1,2\}}(30;2,3)$.}
\label{fig:G-30-2-3-G-30-4-9}
\end{figure}

Let us define the directed graph $\mathcal{G}^{(\mathcal{E})}(n; r, s) := \bigcup_{e \in \mathcal{E}} \mathcal{G}(n; r^e, s^e)$.
For an example, see Figure~\ref{fig:G-30-2-3-G-30-4-9}.

Since $1 \in \mathcal{E}$, it is easy to check that $\mathcal{G}^{(\mathcal{E})}(n; r, s)$ is the graph obtained from $\mathcal{G}(n; r, s)$ by adding a directed edge $v_1 \to v_2$ between each pair $v_1 < v_2$ of vertices of $\mathcal{G}(n; r, s)$ that have distance $e$, for every $e \in \mathcal{E}$.
In particular, this process connects only vertices that are already connected.
Hence, the number of connected components of $\mathcal{G}^{(\mathcal{E})}(n; r, s)$ that have exactly $i$ vertices is equal to the number of connected components of $\mathcal{G}(n; r, s)$ that have exactly $i$ vertices, which is the number $c_i$ that we already determined in the proof of Theorem~\ref{thm:prob}.
Moreover, the probability that a connected component of $\mathcal{G}^{(\mathcal{E})}(n; r, s)$ having exactly $i$ vertices has no adjacent vertices both belonging to $\mathcal{A}$ is equal to the probability $\beta_i^{(\mathcal{E})}$ that the random binary string $\chi_1\cdots\chi_i$ does not contain the substring $10^{e-1}1$, for all $e \in \mathcal{E}$, where $\{\chi_k\}_{k=1}^\infty$ is a sequence of independent identically distributed random variables in $\{0,1\}$ with $\mathbb{P}(\chi_k = 1) = \alpha$.
At this point the same reasonings of~\eqref{equ:chain} yield~\eqref{equ:prob-powers}.
Let us prove that $\big\{\beta_i^{(\mathcal{E})}\big\}_{i=0}^\infty$ is a linear recurrence.

Suppose that $\mathcal{E}$ is finite and let $m := \max(\mathcal{E}) + 1$.
Then $\beta_i^{(\mathcal{E})}$ can be determined by considering a Markov chain.
The states are the binary strings $x_1 \cdots x_m \in \{0,1\}^m$ not containing the substring $10^{e-1}1$, for every $e \in \mathcal{E}$, and one absorbing state. 
A transition from state $x_1 \cdots x_m$ to state $x_2 \cdots x_{m-1} 1$, respectively from state $x_1 \cdots x_m$ to state $x_2 \cdots x_{m-1} 0$, happens with probability $\alpha$, respectively $1 - \alpha$, and all the other transitions are to the absorbing state.
Finally, the probability of $x_1 \cdots x_m$ being the initial state is $\alpha^{x_1 + \cdots + x_m} (1 - \alpha)^{m - (x_1 + \cdots + x_m)}$.
Therefore, letting $u$ be the number of states, we have that $\beta_i^{(\mathcal{E})} = \bm{\pi} \bm{\Sigma}^i (1,1,\dots,1,0)^\text{t}$ for all integers $i \geq 0$, where $\bm{\pi}$ is a row vector of length $u$, $\bm{\Sigma}$ is a $u \times u$ stochastic matrix, and $(1,1,\dots,1,0)^\text{t}$ is column vector of length $t$, assuming the states are ordered so that the absorbing state is the last one.
Consequently, $\big\{\beta_i^{(\mathcal{E})}\big\}_{i=0}^\infty$ is a linear recurrence whose characteristic polynomial is given by the characteristic polynomial of $\bm{\Sigma}$.

Now suppose that $\mathcal{E}$ is cofinite and let $\ell$ be the minimal positive integer such that $e \in \mathcal{E}$ for all integers $e \geq \ell$.
If $\chi_1 \cdots \chi_i$ does not contain $10^{e-1}1$, for every $e \in \mathcal{E}$, then the distance between each pair of $1$s in $\chi_1 \cdots \chi_i$ is less than $\ell$ positions.
In particular, the number of $1$s in $\chi_1 \cdots \chi_i$ is at most $\ell + 1$.
Therefore, for $i \geq \ell + 1$, by elementary probability calculus we can write $\beta_i^{(\mathcal{E})}$ as a linear combination, whose coefficients do not depend on $i$, of the power sums $(i - k + 1) (1 - \alpha)^{i - k}$ where $k=0,\dots,\ell+1$.
Consequently, $\beta_i^{(\mathcal{E})} = (1 - \alpha)^{i - \ell - 1} \Beta^{(\mathcal{E})}(i)$ for some $\Beta^{(\mathcal{E})}(X) \in \mathbb{R}[X]$.
Hence, $\big\{\beta_i^{(\mathcal{E})}\big\}_{i=0}^\infty$ is a linear recurrence and $\gamma_i^{(\mathcal{E})} = \Beta(i - 1) \Beta(i + 1) / \Beta(i)^2$ is a rational function of $i$.

The proof is complete.

\section{General case}\label{sec:generalizations}

As mentioned in the introduction, providing a general formula for the probability of the event $\bigvee_{i=1}^k \big(r_i / s_i \in \mathcal{A} /\! \mathcal{A}\big)$, where $r_1/s_1, \dots, r_k / s_k$ are rational numbers, seems very complex.
In light of the previous reasonings, this task amounts to study the graph $\mathcal{G} := \bigcup_{i = 1}^k \mathcal{G}(n; r_i, s_i)$.
Precisely, one has to classify the connected components of $\mathcal{G}$, and to determine the probability that each of them does not have two adjacent vertices both belonging to $\mathcal{A}$.

\begin{figure}[h]
\includegraphics[width=\textwidth]{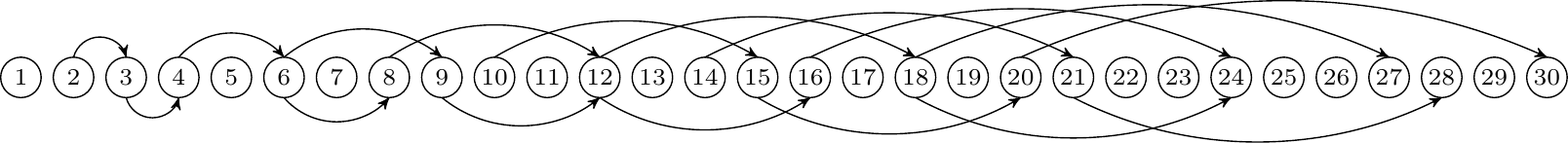}
\caption{The directed graph $\mathcal{G}(30;2,3) \cup \mathcal{G}(30;3,4)$.}
\label{fig:G-30-2-3-G-30-3-4}
\end{figure}

\begin{figure}[h]
\includegraphics[width=0.538\textwidth]{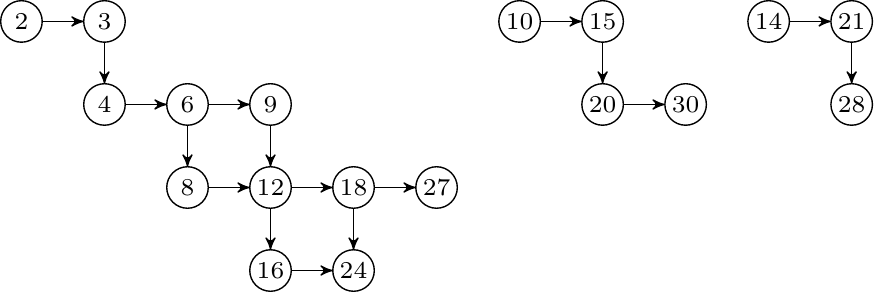}
\caption{The connected components of $\mathcal{G}(30;2,3) \cup \mathcal{G}(30;3,4)$ that have at least $2$ vertices.
Each horizontal, respectively vertical, edge corresponds to multiply the value of a vertex by $3/2$, respectively $4/3$.}
\label{fig:G-30-2-3-G-30-3-4-comp}
\end{figure}

\begin{figure}[h]
\includegraphics[width=\textwidth]{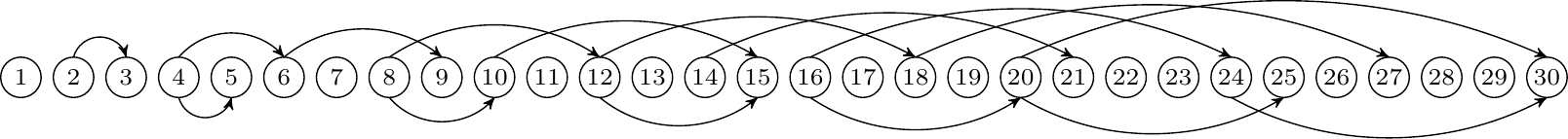}
\caption{The directed graph $\mathcal{G}(30;2,3) \cup \mathcal{G}(30;4,5)$.}
\label{fig:G-30-2-3-G-30-4-5}
\end{figure}

\begin{figure}[h]
\includegraphics[width=0.64\textwidth]{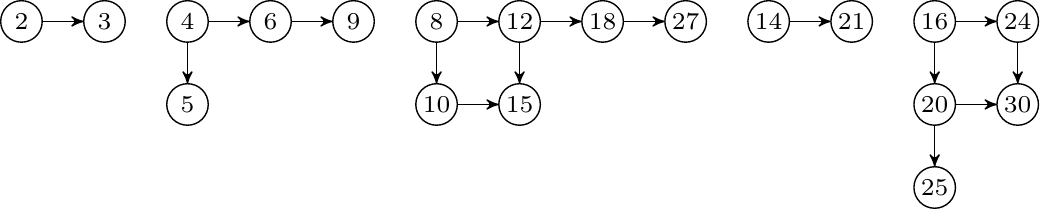}
\caption{The connected components of $\mathcal{G}(30;2,3) \cup \mathcal{G}(30;4,5)$ that have at least $2$ vertices.
Each horizontal, respectively vertical, edge corresponds to multiply the value of a vertex by $3/2$, respectively $5/4$.}
\label{fig:G-30-2-3-G-30-4-5-comp}
\end{figure}

\FloatBarrier 

If the multiplicative group generated by $\{r_i / s_i\}_{i=1}^k$ is cyclic, then the connected components of $\mathcal{G}$ have a somehow ``linear'' structure, and proving formulas similar to~\eqref{equ:prob} and \eqref{equ:prob-powers} is doable.

If the generated group has rank $R > 1$, then each connected component of $\mathcal{G}$ is isomorphic to a subgraph of the $R$-dimensional grid graph.
For examples, see Figures~\ref{fig:G-30-2-3-G-30-3-4}, \ref{fig:G-30-2-3-G-30-3-4-comp}, \ref{fig:G-30-2-3-G-30-4-5}, and~\ref{fig:G-30-2-3-G-30-4-5-comp}.

\section{Visible lattice points}

Another direction of research can be generalizing ratio sets to sets of visible lattice points.
Let $d \geq 2$ be an integer.
For every $\mathcal{A} \subseteq \mathbb{N}$, a lattice point $P \in \mathbb{N}^d$ is said to be visible in the lattice $\mathcal{A}^d$ if the line segment from $\mathbf{0} \in \mathbb{Z}^d$ to $P$ intersects $\mathcal{A}^d$ only in $P$.
Let $\vis(\mathcal{A}^d)$ be the set of lattice points visible in $\mathcal{A}^d$.
There is a natural bijection between $\vis(\mathcal{A}^2)$ and $\mathcal{A} /\! \mathcal{A}$, given by $(x_1, x_2) \mapsto x_1 / x_2$.
Hence, $\vis(\mathcal{A}^d)$ can be considered as a $d$-dimensional generalization of the ratio set $\mathcal{A} /\! \mathcal{A}$ (see also~\cite{MR4092824} for a similar generalization of ratio sets).

\begin{figure}[h]
\includegraphics[width=0.93\textwidth]{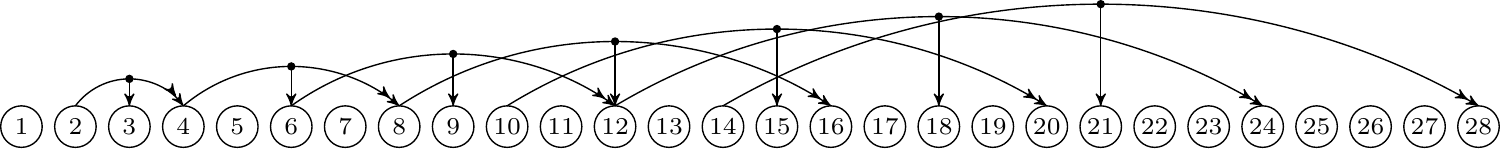}
\caption{The hypergraph $\mathcal{H}(28;2,3,4)$.}
\label{fig:hypergraph}
\end{figure}

Cilleruelo and Guijarro-Ord\'{o}\~{n}ez~\cite{MR3649012} gave an asymptotic formula for the cardinality of $\vis(\mathcal{A}^d)$ for $\mathcal{A} \in \mathcal{B}(n, \alpha)$.
A natural question is if Theorem~\ref{thm:prob} can be generalized to a formula for $\mathbb{P}\big((x_1, \dots, x_d) \in \vis(\mathcal{A}^d)\big)$, where $(x_1, \dots, x_d) \in \mathbb{N}^d$.
This amount to study the hypergraph $\mathcal{H}(n; x_1, \dots, x_d)$ defined as having vertices $1, \dots, n$ and hyperedges $(x_1 t, \dots, x_d t)$, for every positive integer $t \leq n /\!\max(x_1, \dots, x_d)$.
For an example, see Figure~\ref{fig:hypergraph}.

\section{Acknowledgements}

The author thanks Paolo Leonetti and Daniele Mastrostefano for suggestions that improved the quality of the article.

\FloatBarrier
\bibliographystyle{amsplain}

\end{document}